\UseRawInputEncoding
\documentclass[11pt,oneside,reqno]{amsart}
\newtheorem{theorem}{Theorem}[section]
\usepackage[left=2.cm,top=1.5cm,bottom=1.5cm,right=2.0cm]{geometry}
\usepackage{latexsym}
\usepackage{amssymb,amsmath,amsfonts}
\usepackage[pagewise]{lineno}
\usepackage{multirow}
\usepackage{graphicx}
\usepackage{cite}
\usepackage[colorlinks = true]{hyperref}
\theoremstyle{definition}
\newtheorem{definition}[theorem]{Definition}

\newtheorem{proposition}[theorem]{Proposition}

\newcommand{\A}{\mathbb A}
\newcommand{\R}{\mathbb R}

\newcommand{\N}{\mathbb N}

\newcommand{\Scal}{\mathcal{S}}

 \global\parskip 6pt \oddsidemargin
0.1cm \evensidemargin 0.4cm \headheight 0.1cm \topmargin -2.0cm
\begin{document}
\pagenumbering{arabic}
\title[Solow-Swan model with nonlocal fractional derivative operator]{Analysis of Solow-Swan model with nonlocal fractional derivative operator}

\author{ Mathew O. Aibinu$^{1,2,3,6*}$, Kevin J. Duffy$^{4}$, Sibusiso Moyo$^{5,6}$}
\address{$^{1}$ Department of Mathematics and Statistics, University of Regina, Regina, SK S4S 0A2, Canada}
\address{$^{2}$ Ingram School of Engineering, Texas State University, TX 78666, United States}
\address{$^{3}$ School of Computer Science and Applied Mathematics, University of the Witwatersrand, Johannesburg 2050, South Africa}
\address{$^{4}$ School of Mathematics, Statistics \& Computer Science, University of KwaZulu-Natal, Durban 4000, South Africa}
\address{$^{5}$Department of Applied Mathematics and School for Data Science and Computational Thinking, Stellenbosch University, South Africa}
\address{$^{6}$National Institute for Theoretical and Computational Sciences (NITheCS), South Africa}
\email{$^*$moaibinu@yahoo.com (MO Aibinu)} 

\begin{abstract} 
The Solow-Swan equation is a foundational model in the evolution of modern economic growth theory. It offers key insights into the long-term behaviour of capital accumulation and output. Since its inception, the model has served as a cornerstone for understanding macroeconomic dynamics and has inspired a vast body of subsequent research. However, traditional formulations of the Solow-Swan model rely on integer-order derivatives, which may not fully capture the memory and hereditary properties often observed in real-world economic systems. In this paper, we extend the classical Solow-Swan framework by incorporating memory effects through the use of fractional calculus. The fractional model accounts for the influence of past states on the present rate of capital change, a feature not accommodated in the standard model. We present a comparative analysis of the capital dynamics under both the classical and fractional-order  formulations of the Solow-Swan equation. 
\end{abstract}


\maketitle

\section{Introduction}\label{sse6}
\par The fundamental differential equation
\begin{equation}\label{sse1}
\frac{dk(t)}{dt}=k(t)\left(pk^{\mu-1}(t)-q\right), \ \ k(0)=k_0,
\end{equation}
is referred to as the Solow-Swan model for capital accumulation that describes the capital ($K$) per labour ($L$)  with ratio, $k = K/L,$ as a function of time, $t,$ where $\mu$ is the exponent and the scaling parameters are $p$ and $q.$ The Solow-Swan model appears as a consequence of modern growth theory in economics and remains an  active research topics (see, e.g, \cite{Brunner1}). It is also present in most undergraduate curricula for economics (see, e.g, \cite{RaghavendraP}). The Solow-Swan model is a dynamic model that has been adapted to consider modelling of exponentially increasing  phenomena of physical systems in areas such as climate change (see, e.g, \cite{BurdaZS}), corruption (see, e.g, \cite{GonzalezCAD, SpyromitrosP}), education (see, e.g, \cite{VikiaYM}), and over-exploitation of natural resources (see, e.g, \cite{VikiaYM}). The Solow-Swan model has been studied extensively due to its wide applications. Brunner et al. \cite{Brunner1} applied the method of least squares to determine the best-fit parameters of the Solow-Swan model (\ref{sse1}). The stability of the spatially homogeneous equilibrium in the Solow–Swan model has been analyzed in various studies (see, e.g, \cite{NetoCJ}). The classical Solow-Swan model predicts convergence to a steady state. In contrast, the endogenous growth-cycle model allows for persistent growth fluctuations or cycles. These cycles arise from feedback loops within the economy. As a result, the model captures more complex and realistic long-term dynamics of capital accumulation and output (see, e.g, \cite{BroitmanC}). A generalization of the Solow-Swan model through the introduction of a non-constant labour growth rate allows for more realistic modelling of economic dynamics by relaxing the assumption of a fixed population or labour force growth. In this extended framework, the labour force may grow at variable rates over time, influenced by demographic trends such as migration, education policies, or labour market conditions. This modification enables the model to capture transitional dynamics more accurately and to reflect real-world fluctuations in labour supply, ultimately leading to more nuanced predictions of long-term growth and convergence behaviour \cite{TolmachevLMB}. Unlike the classical Solow-Swan model with ordinary integer-order derivatives that has received much attention, the literature on the Solow-Swan model with fractional derivatives is still regarded as rare (see, e.g, \cite{Traore}). Fractional Differential Equations (FDEs)  are often used to model growth associated with memory effects. The memory effects are due to the non-local properties that are possessed by the FDEs. This is an edge that fractional derivatives have over classical derivatives, where the effect is generally ignored. The FDE is a suitable concept for modelling the growth of many economical processes because many economical processes have memory effect in their nature (see, for example, \cite{Johansyah1}).

\par The motivation for this study stems from the fact that most existing research on the Solow-Swan equation has focused on models involving ordinary integer-order derivatives. In contrast, formulations of the Solow-Swan equation that incorporate memory effects have received relatively little attention. This paper aims to address that gap by introducing memory effects into the model through the use of a fractional-order derivative. A comparative analysis is conducted to examine the behaviour of capital dynamics under both the classical integer-order and fractional-order formulations of the Solow-Swan equation. The study investigates how the inclusion of a fractional derivative influences capital behaviour. Additionally, the effects of key scaling parameters on the dynamics of the models are explored.

\section{Preliminaries}
Some significant definitions and lemmas that are essential to establish the results in this paper are presented in this section.
\begin{definition}
Over the set
$$\A=\left\{k(t): \exists \ Q, \eta_1, \eta_2 >0, |k(t)|<Q\exp{\left({|t|}/\eta_j\right)}, \mbox{if} \ t\in (-1)^j\times [0, \infty) \right\},$$
the Sumudu Transform (ST) is defined as (see, e.g, \cite{Belgacem2})
\begin{equation}\label{em1}
\Scal [k(t)]=\int^{\infty}_{0}k(tu)e^{-t}dt, \ u\in (-\eta_1, \eta_2).
\end{equation}
The ST satisfies the linearity conditions (see, e.g, \cite{Moltot1, Belgacem2, Watugala1, Belgacem1}) and it is a method that is highly cherished for its unit preserving and domains scaling properties \cite{TunG1}. Denoting $S\left[k(t)\right]$ by $K(u),$ the ST for the $n^{th}$-order derivative is
\begin{equation}\label{sumud20}
S\left[k^n(t)\right]=\frac{1}{u^n}\left[K(u)-\displaystyle\sum_{i=0}^{n-1}u^ik^i(t)|_{t=0}\right].
 \end{equation}
Therefore, for the first order derivative, it is simply
 \begin{equation}\label{sumud19}
S\left[k'(t)\right]=\frac{1}{u}\left[K(u)-k(0)\right].
 \end{equation}
\end{definition}
\begin{definition}
 Let $a>0, b>0$ be positive real numbers. The Caputo derivative of order $\alpha$ is defined as
$$^C_aD^{\alpha}p(t)=\frac{1}{\Gamma\left(1-\alpha\right)}\int^{t}_{a}(t-\eta)^{-\alpha}p'(\eta)d\eta,$$
where $0<\alpha <1.$ The Caputo derivative admits the ST in the form (see, e.g, \cite{Bodkhe})
\begin{equation}\label{em5}
\Scal \left[^C_0D^{\alpha}k(t)\right]=u^{-\alpha}\left(K(u)-k(0)\right).
\end{equation}
\end{definition}

\begin{proposition}
Let $\psi, \varphi :[0, \infty)\rightarrow \R,$ the classical convolution product is given by
$$(\psi * \zeta)(t)=\int^{t}_{0}\psi(t-x)\zeta(x)dx.$$
The ST for the convolution product is given by
\begin{eqnarray*}
\Scal \left[(\psi * \zeta)(t)\right]&=&u\Scal[\psi(t)]\Scal[\zeta(t)]\\
&=&u\psi(u)\zeta(u).
\end{eqnarray*}
\end{proposition}

\begin{definition}\label{em20}
The Mittag-Leffler function $E_{\mu}(t)$ is defined as
$$E_{\alpha}(t)=\sum_{n=0}^{\infty}\frac{t^n}{\Gamma\left(\alpha n+1\right)}, \alpha >0.$$
The following results about Mittag-Leffler functions and ST are well known (see, e.g, \cite{Nanware}):
\begin{itemize}
	\item [(i)]$\Scal\left[E_{\alpha}\left(-at^{\alpha}\right)\right]=\frac{1}{1+au^{\alpha}},$
	\item [(ii)]$\Scal\left[1-E_{\alpha}\left(-at^{\alpha}\right)\right]=\frac{au^{\alpha}}{1+au^{\alpha}}.$
\end{itemize}
\end{definition}

\begin{definition}
The power series method expresses the solution of a differential equation as an infinite sum. It is an effective approach, especially applicable when standard elementary methods are insufficient. Consider an equation whose solution is $k(t)$ and that contains a nonlinear term $N[k].$ Suppose the solution is decomposed as 
$$k(t)=\sum_{n=0}^{\infty}k_nt^n$$ and the nonlinear term is expressed as 
$$N[k(t)]=\sum_{n=0}^{\infty}A_nt^n,$$
where $A_n$ are special polynomials that are being referred to as the Adomian polynomials \cite{Adomian1}. The $A_n$ polynomials are defined by \cite{Adomian2}
$$A_n=\frac{1}{n!}\left[\frac{d^n}{d{x}^n}f\left(\displaystyle\sum_{i=0}^{\infty}{x}^ik_i\right)\right]\bigg|_{x=0}.$$
The first few terms of the Adomian polynomials are generated as
$$\begin{cases}
A_0=f(k_0),\\
A_1=k_1f'(k_0),\\
A_2=k_2f'(k_0)+\frac{k_1^2}{2!}f''(k_0)\\
A_3=k_3f'(k_0)+k_1k_2f''(k_0)+\frac{k_1^3}{3!}f'''(k_0)\\
\vdots
\end{cases}$$
Observe that the polynomials $A_n,$ are generated for each nonlinearity such that $A_0$ depends only on $k_0,$ $A_1$ depends only on $k_0,$ and $k_1,$ $A_2$ depends on $k_0,$ $k_1,$ $k_2,$ etc.
\end{definition}
\section{Main results}
The Solow-Swan equation (\ref{sse1}) is nonlinear in nature, and such equations often resist well-known analytical methods. In this study, the solutions of the Solow–Swan equation are considered for both the integer-order and Caputo fractional-order cases. The objective is to obtain approximate analytical solutions using a hybrid of the ST method (see, e.g, \cite{LiuC, Aibinu8, Aibinu7, Aibinu6}). The Solow–Swan model with a memory effect is also examined, and the impact of the fractional-order derivative on the dynamics of capital is analyzed.

\subsection{The Solow-Swan model with ordinary integer-order derivative}
\par
\begin{figure}
\centering
\includegraphics[width=9.0cm ,height=6.0cm]{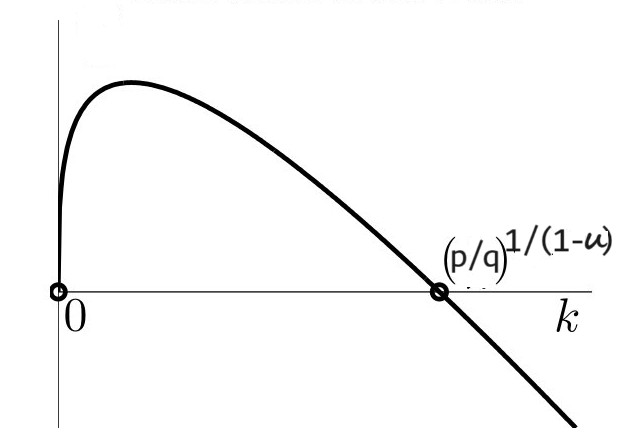}
 \caption{The graph showing the points of equilibria of (\ref{sse1}).}
 \label{sse12}
\end{figure}

 Observe that (\ref{sse1}) is a nonlinear equation. At ${dk}/{dt}=0,$ the points of equilibria of (\ref{sse1}) occur at $k=0, \left(p/q\right)^\frac{1}{1-\mu}$ (see  Figure \ref{sse12} ). Changing the value of either $p$ or $q$ in Figure \ref{sse12} alters the scale and the numerical value of the non-zero equilibrium; however, the qualitative shape of the graph remains invariant. The ST of (\ref{sse1}) is taken as
$$\Scal \left[dk/dt\right]=p\Scal \left[k^{\mu}\right]-q\Scal \left[k\right],$$
and applying the relation (\ref{em5}), it leads to
$$K(u)-k(0)=u\left(p\Scal \left[k^{\mu}\right]-q\Scal \left[k\right]\right),$$
where $\Scal [k]$ is denoted by $K(u).$ The variational iteration formula takes the form
\begin{eqnarray}\label{sse7}
K_{n+1}(u)&=&K_n(u)+\lambda(u)\bigg(\frac{K_n(u)-k_0}{u}-p\Scal \left[N[k_n]\right] + q\Scal \left[k_n\right]\bigg), n\in \N,
\end{eqnarray}
where $N[k_n]= k_n^{\mu}.$ Treating $q\Scal \left[k_n\right]-p\Scal \left[N[k_n]\right]$ as the restricted term in (\ref{sse7}) and taking its classic variation operators gives the Langrange multiplier as 
\begin{equation}\label{sse8}
\lambda(u)=-u.
\end{equation}
Substituting the expression (\ref{sse8}) into (\ref{sse7}) and taking its inverse ST gives
$$k_{n+1}(t)=k_0+\Scal^{-1}\left[u\left(p\Scal \left[N[k_n]\right] - q\Scal \left[k_n\right]\right)\right],$$
where $k_0(t)=k_0.$ Let 
\begin{equation}\label{sse9}
k_n=\sum_{i=0}^nw_i, 
\end{equation}
and the decomposition of the nonlinear term becomes 
\begin{equation}\label{sse10}
N[k_n]=\sum_{i=0}^{n}A_i, \ \mbox{with} \ A_i=\frac{1}{i!}\left[\frac{d^i}{d{x}^i}f\left(\sum_{n=0}^{\infty}{x}^nw_n\right)\right]\bigg|_{x=0},
\end{equation}
where $A_i$ are the Adomian polynomials. Then the Adomian series of the $w^{\alpha}$ is as follows:
\begin{equation}\label{sse11}
\begin{cases}
A_0=w_0^{\mu},\\
A_1={\mu}w_1w_0^{\mu-1},\\
A_2={\mu}w_2w_0^{\mu-1}+\mu(\mu-1)\frac{w_1^2}{2!}w_0^{\mu-2},\\
\vdots
\end{cases}
\end{equation}\

\begin{figure}
\centering
\includegraphics[width=9.0cm ,height=5.5cm]{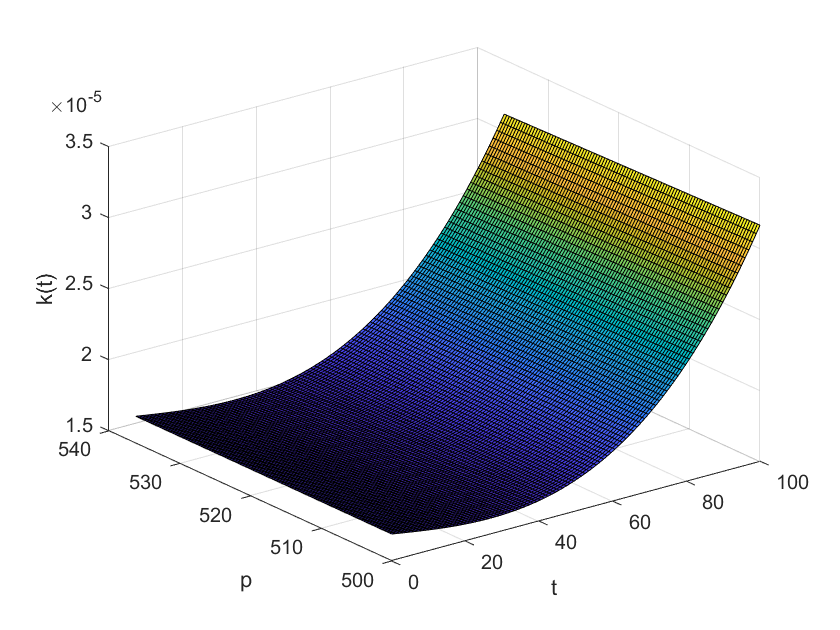}
 \caption{Graph of (\ref{sse13}) when $q$ and $\mu$ are held constant.}
 \label{sse15}
\end{figure}

\begin{figure}
\centering
\includegraphics[width=9.0cm ,height=5.5cm]{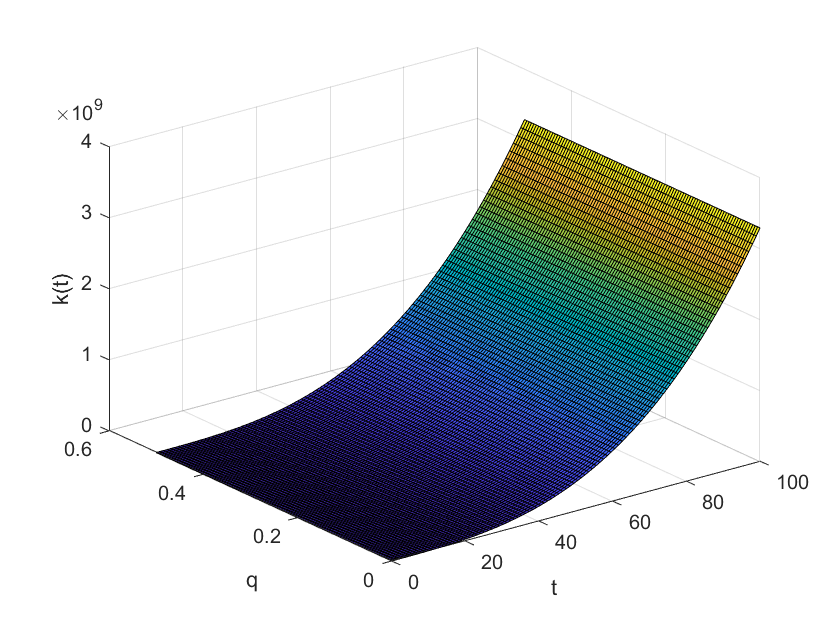}
\caption{Graph of (\ref{sse13}) when $p$ and $\mu$ are held constant.}
 \label{sse16}
\end{figure}

\begin{figure}
\centering
\includegraphics[width=9.0cm ,height=5.5cm]{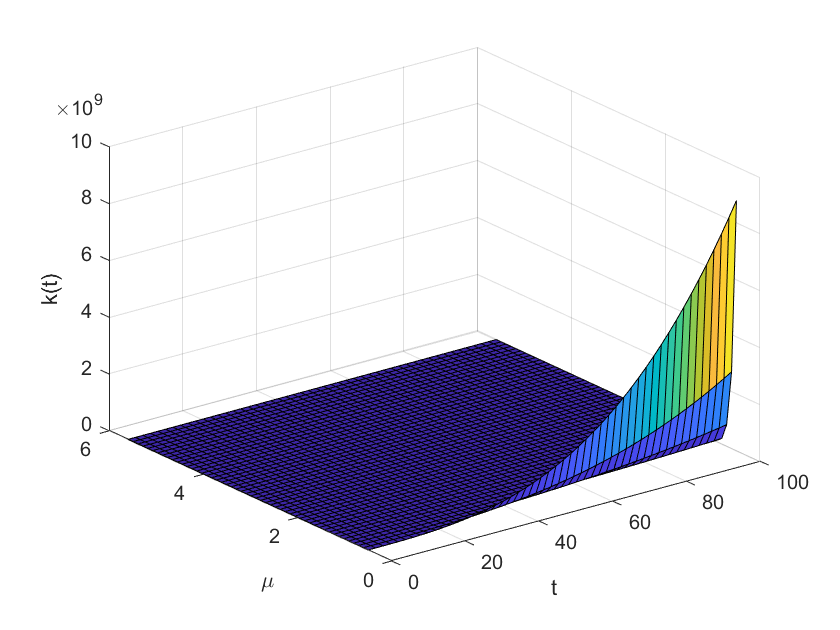}
\caption{Graph of (\ref{sse13}) when $p$ and $q$ are held constant.}
 \label{sse17}
\end{figure}

Therefore, the iteration formula is derived as follows
$$\begin{cases}
w_0(t)=k(0)=k_0,\\
w_{n+1}(t)= \Scal^{-1}\left[u\left(p\Scal \left[A_n\right] - q\Scal \left[w_n\right]\right)\right],
\end{cases}$$
which produces the sequence
$$\begin{cases}
w_0&=k_0,\\
w_1&=\left(pk_0^{\mu}-qk_0\right)t,\\
w_2&=\left(pk_0^{\mu}-qk_0\right)\left(p\mu k_0^{\mu-1}-q\right)\frac{t^{2}}{2!},\\
w_3&=\left(pk_0^{\mu}-qk_0\right)\left\{\left(p\mu k_0^{\mu-1}-q\right)^2-\frac{p\mu(\mu-1)}{2}k_0^{\mu-2} \right\}\frac{t^{3}}{3!},\\
\vdots
\end{cases}$$
Therefore, the solution is given by
\begin{eqnarray}\label{sse13}
k(t)&=&\displaystyle \lim_{n\rightarrow \infty}k_n=\displaystyle \lim_{n\rightarrow \infty}\sum_{n=0}^{\infty}w_n=k_0+\left(pk_0^{\mu}-qk_0\right)t\nonumber\\
&&+\left(pk_0^{\mu}-qk_0\right)\left(p\mu k_0^{\mu-1}-q\right)\frac{t^{2}}{2!}\\
&&+\left(pk_0^{\mu}-qk_0\right)\left\{\left(p\mu k_0^{\mu-1}-q\right)^2-\frac{p\mu(\mu-1)}{2}k_0^{\mu-2} \right\}\frac{t^{3}}{3!}+...\nonumber
\end{eqnarray}
Using the statistical data from the literature \cite{Brunner1}, we assigned suitable values to the parameters in \ref{sse1} to display the graphs of its solution. In Figure \ref{sse15}, equation (\ref{sse13}) is plotted with $q$ held constant as $p$ varies. Similarly, in Figure \ref{sse16}, $p$ is held constant while $q$ varies. Figures \ref{sse15} and  \ref{sse16} show that the qualitative shape of equation (\ref{sse13}) remains invariant when the values of either $p$ or $q$ are changed. The observed changes occur only in the scale, corresponding to variations in the parameters $p$ and $q.$ Figure \ref{sse17} displays the graph of equation (\ref{sse13}) for values of $\mu$ in the range $[0.5, 5.6].$ 

\subsection{The Solow-Swan model with Caputo derivative} The Caputo fractional derivative is a widely used form of fractional differentiation, notable for naturally accommodates initial conditions in a manner similar to integer-order derivatives (see, e.g, \cite{OdibatZ}). This makes it particularly useful in modelling problems where non-local properties and historical interactions must be taken into account. Under some invariability assumptions on the function, it shares dual relationships with the Riemann–Liouville fractional derivative. The Solow-Swan model in the sense of Caputo is given as
\begin{equation}\label{sse2}
^C_0D^{\alpha}k(t)=k(t)\left(pk^{\mu-1}(t)-q\right), \ \ k(0)=k_0.
\end{equation}


Observe that when $\alpha=1,$ the Caputo derivative converges to the classical ordinary integer-order derivative. The ST of (\ref{sse2}) is taken as
$$\Scal \left[^C_0D^{\alpha}k\right]=p\Scal \left[k^{\mu}\right]-q\Scal \left[k\right].$$
Applying the relation (\ref{em5}) leads to
$$\frac{\Scal [k]-k_0}{u^{\alpha}}=p\Scal \left[k^{\mu}\right]-q\Scal \left[k\right].$$
The variational iteration formula takes the form
\begin{eqnarray}\label{sse3}
K_{n+1}(u)&=&K_n(u)+\lambda(u)\bigg(\frac{K_n(u)-k_0}{u^{\alpha}}-p\Scal \left[N[k_n]\right] + q\Scal \left[k_n\right]\bigg), n\in \N,
\end{eqnarray}
where $N[k_n]= k_n^{\mu}.$ Treating $q\Scal \left[k_n\right]-p\Scal \left[N[k_n]\right]$ as the restricted term in (\ref{sse3}) and taking its classic variation operators gives the Langrange multiplier as 
\begin{equation}\label{sse4}
\lambda(u)=-u^{\alpha}.
\end{equation}
Substituting the expression (\ref{sse4}) into (\ref{sse3}) and taking its inverse Sumudu transform gives
$$k_{n+1}(t)=k_0+\Scal^{-1}\left[u^{\alpha}\left(p\Scal \left[N[k_n]\right] - q\Scal \left[k_n\right]\right)\right],$$
where $k_0(t)=k_0.$ Let $k_n, N[k_n]$ and the Adomian series be given by (\ref{sse9}), (\ref{sse10}) and (\ref{sse11}) respectively. Therefore, the iteration formula is derived as follows
$$\begin{cases}
w_0(t)=k(0)=k_0,\\
w_{n+1}(t)= \Scal^{-1}\left[u^{\alpha}\left(p\Scal \left[A_n\right] - q\Scal \left[w_n\right]\right)\right],
\end{cases}$$
which produces the sequence
$$\begin{cases}
w_0&=k_0,\\
w_1&=\left(pk_0^{\mu}-qk_0\right)\frac{t^{\alpha}}{\Gamma(\alpha+1)},\\
w_2&=\left(pk_0^{\mu}-qk_0\right)\left(p\mu k_0^{\mu-1}-q\right)\frac{t^{2 \alpha}}{\Gamma(2 \alpha+1)},\\
w_3&=\left(pk_0^{\mu}-qk_0\right)\left\{\left(p\mu k_0^{\mu-1}-q\right)^2-\frac{p\mu(\mu-1)}{2}k_0^{\mu-2} \right\}\frac{t^{3\alpha}}{\Gamma(3\alpha +1)},\\
\vdots
\end{cases}$$

Therefore, the solution is given by
\begin{equation}\label{sse18}
\begin{split}
k(t)&=\displaystyle \lim_{n\rightarrow \infty}k_n=\displaystyle \lim_{n\rightarrow \infty}\sum_{n=0}^{\infty}w_n=k_0+\left(pk_0^{\mu}-qk_0\right)\frac{t^{\alpha}}{\Gamma(\alpha+1)}\\
&+\left(pk_0^{\mu}-qk_0\right)\left(p\mu k_0^{\mu-1}-q\right)\frac{t^{2 \alpha}}{\Gamma(2 \alpha+1)}\\
&+\left(pk_0^{\mu}-qk_0\right)\left\{\left(p\mu k_0^{\mu-1}-q\right)^2-\frac{p\mu(\mu-1)}{2}k_0^{\mu-2} \right\}\frac{t^{3\alpha}}{\Gamma(3\alpha +1)}+...
\end{split}
\end{equation}

\begin{figure}
\centering
\includegraphics[width=9.0cm ,height=5.5cm]{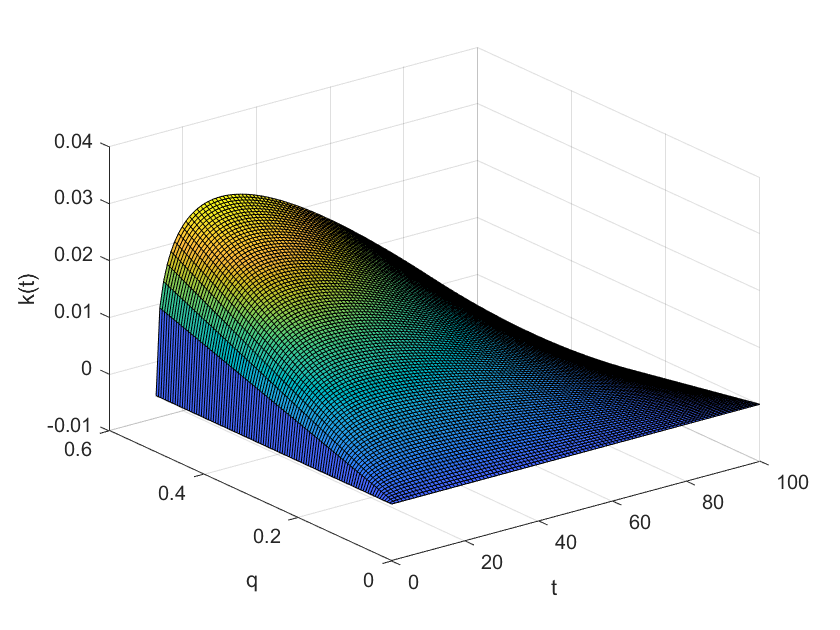}
 \caption{Graph of (\ref{sse18})  for constant values of $\alpha, \mu$ and $p.$}
 \label{sse19}
\end{figure}

\begin{figure}
\centering
\includegraphics[width=9.0cm ,height=5.5cm]{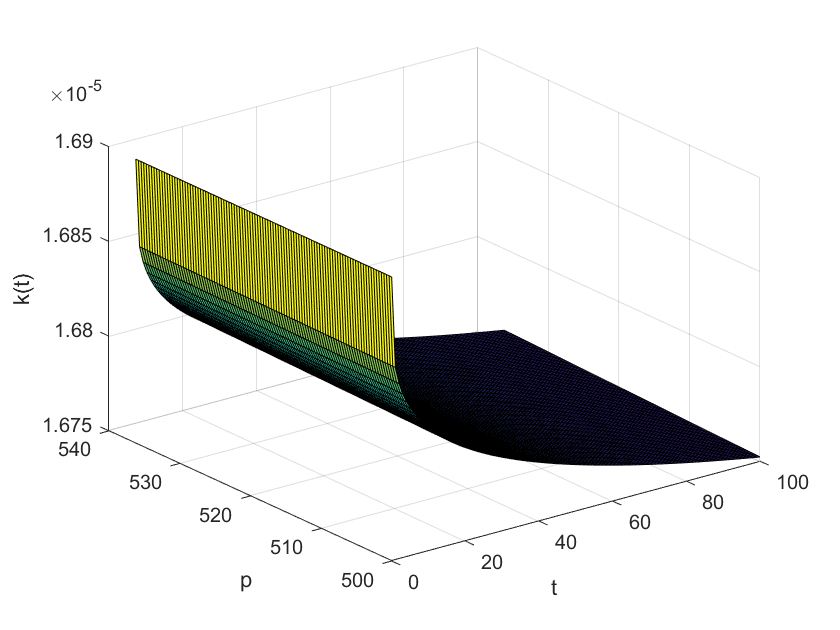}
\caption{Graph of (\ref{sse18}) for constant values of $\alpha, \mu$ and $q.$}
 \label{sse20}
\end{figure}

\begin{figure}
\centering
\includegraphics[width=9.0cm ,height=5.5cm]{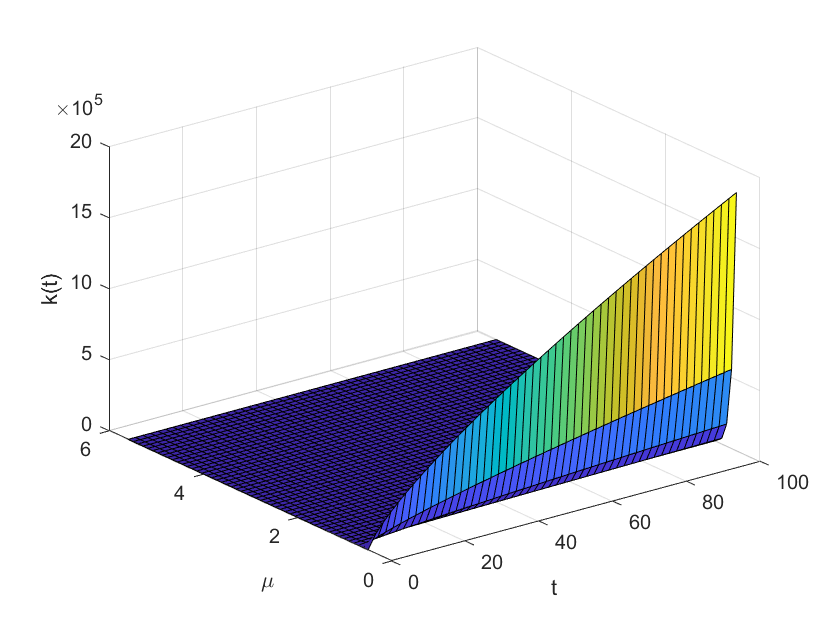}
\caption{Graph of (\ref{sse18}) for constant values of $\alpha, p$ and $q.$}
 \label{sse21}
\end{figure}

\begin{figure}
\centering
\includegraphics[width=9.0cm ,height=5.5cm]{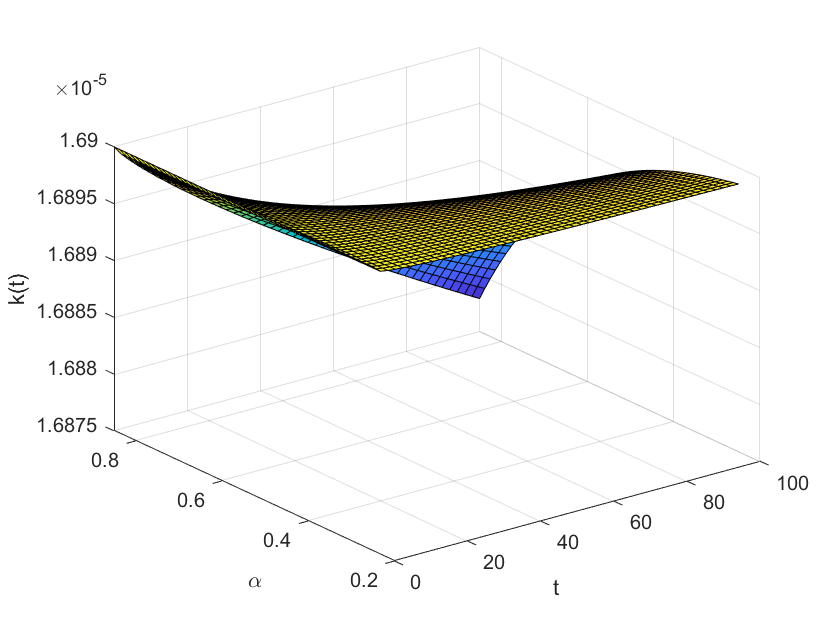}
\caption{Graph of (\ref{sse18}) for constant values of $\mu, p$ and $q.$}
 \label{sse22}
\end{figure}

The Caputo derivative in (\ref{sse2}) provides more degrees of freedom in the experimental simulations and allows for the  consideration of  memory effects. The solution of (\ref{sse2}) is given by (\ref{sse18}).  In Figure \ref{sse19}, (\ref{sse18}) is plotted with $\alpha, \mu$ and $p$ held constant while $q$ varies. Similarly, Figure \ref{sse20} illustrates equation (\ref{sse18}) with $\alpha, \mu$ and $q$ held constant while $p$ varies. Figure \ref{sse21} illustrates the graph of (\ref{sse18}) with $\alpha, p$ and $q$ held constant, and $\mu$ varying. Figures \ref{sse19}, \ref{sse20} and \ref{sse21} depict the Solow-Swan model in the context of a fractional derivative of order $\alpha=0.25.$ Figure \ref{sse22} on the other hand, shows the model with a fractional derivative order ranging from $\alpha\in[0.3, 0.9],$ with $p, q$ and $\mu$ held constant. The Caputo derivative incorporates the effects of past states and captures the non-local characteristics inherent in the model.

\subsection{Analysis in terms of the capital $K$ and the labour $L$}
The equilibrium point for equations (\ref{sse1}) and (\ref{sse2}) is $k = 0, \left(p/q\right)^\frac{1}{1-\mu}.$ The equilibrium $k = 0$ is unstable because for a small increase at the point, that is $k>0,$ consequently ${dk}/{dt}>0,$ so $k$ will increase. At $k=\left(p/q\right)^\frac{1}{1-\mu},$ the graphs of $k(t)$ will have an inflexion point. In other words, the right-hand sides of (\ref{sse1}) and  (\ref{sse2}) will reach their respective maximum values. Equations (\ref{sse1}) and (\ref{sse2}) exhibit asymptotic stability at $k=\left(p/q\right)^\frac{1}{1-\mu}$ with $k(t)$ converging to the non-zero equilibrium as $t$ approaches infinity. Recall that the Solow-Swan models for capital accumulation describe the capital ($K$) per labour ($L$) with ratio, $k = K/L,$ where $L(t)=L_0e^{\psi t}.$ Whenever $k(t)$ converges asymptotically to a stable equilibrium $k_1,$ then $K(t)$ inevitably asymptotically converges to $k_1L(t).$ The models predict that, in the long term, both capital and effective labour grow exponentially and in tandem. This indicates that under steady-state conditions and constant returns to scale, the economy approaches a balanced growth path. In this state, output, capital, and effective labour expand at the same exponential rate, and key ratios remain constant. If the initial capital stock is significantly below its steady-state level, the model predicts a period of rapid capital accumulation. When capital is scarce relative to labour, the marginal productivity of capital is high, which actuates investment and accelerates capital growth. Over time, as capital increases and approaches a level proportional to effective labour, its growth rate stabilizes, and the economy converges to its long-run balanced path. Eventually, the system settles into a long-run trajectory where capital and labour maintain a constant ratio, growing proportionally over time. This convergence to proportional growth reflects the models’ inherent tendency toward equilibrium. This long-term behaviour not only illustrates the self-correcting nature of the Solow-Swan model but also highlights the importance of initial conditions in determining short-run dynamics. Ultimately, the model provides a powerful framework for understanding how economies adjust over time and how balanced growth emerges from the interaction between capital accumulation, labour expansion, and productivity.

\vspace{1.0cm}
\subsection{Conclusion} The Caputo derivative is widely regarded as the most suitable fractional operator for modelling real-world phenomena. It accommodates traditional initial and boundary conditions in a manner consistent with classical differential equations. As a result, it has become a preferred tool among researchers for formulating and analyzing problems that exhibit memory and hereditary properties. In this paper, we examine the Solow-Swan equation, a representative nonlinear model in economic growth theory that often challenges conventional analytical methods. To enhance the accuracy of the model, memory effects are introduced by reformulating the equation using the Caputo fractional derivative. This modification enables the model to account for both the influence of historical states and non-local behaviour. The new model offers a more comprehensive representation of capital dynamics over time.

A hybrid of the ST method is employed to derive solutions to the Solow-Swan equation in both its classical and fractional  forms. Numerical simulations are carried out using MATLAB to visualize and compare the behaviours of the two models. The resulting graphical representations clearly illustrate the impact of incorporating memory effects through fractional calculus into the dynamics of economic growth. The models predict that, in the long term, both capital and labour grow exponentially and simultaneously. This outcome aligns with the concept of a balanced growth path, where the key inputs to production, which are capital and labour, expand at steady exponential rates. In the traditional Solow-Swan framework, exponential growth in labour is typically driven by population growth, while capital accumulation depends on the savings rate, depreciation, and labour force growth. When both capital and labour grow exponentially, the economy tends to move along a stable trajectory. Exponential growth in both capital and labour suggests that the economy is operating on a stable and predictable growth path. This has practical implications for long-term investment decisions, fiscal planning, and infrastructure development.

\end{document}